\documentclass[12pt]{article} \usepackage{amsfonts} \usepackage{amssymb}\usepackage{amsmath}
\makeatletter
\renewcommand{\@makefnmark}{}
\makeatother
\begin{document}
\baselineskip=14pt
\pagestyle{plain}
{\Large
\newcommand{\mun}{\mu_{n,j}}
\newcommand{\wnj}{W_{N,j}(\mu)}
\newcommand{\ntm}{|2n-\theta-\mu|}
\newcommand{\pwp}{PW_\pi^-}
\newcommand{\cpm}{\cos\pi\mu}
\newcommand{\spm}{\frac{\sin\pi\mu}{\mu}}
\newcommand{\tet}{(-1)^{\theta+1}}
\newcommand{\dm}{\Delta(\mu)}
\newcommand{\dl}{\Delta(\lambda)}

\newcommand{\dnm}{\Delta_N(\mu)}
\newcommand{\sni}{\sum_{n=N+1}^\infty}
\newcommand{\muk}{\sqrt{\mu^2+q_0}}
\newcommand{\agt}{\alpha, \gamma, \theta,}
\newcommand{\muq}{\sqrt{\mu^2-q_0}}
\newcommand{\smn}{\sum_{n=1}^\infty}
\newcommand{\lop}{{L_2(0,2\pi)}}
\newcommand{\tdm}{\tilde\Delta_+(\mu_n)}
\newcommand{\dmn}{\Delta_+(\mu_n)}
\newcommand{\dd}{D_+(\mu_n)}
\newcommand{\ddd}{\tilde D_+(\mu_n)}
\newcommand{\sn}{\sum_{n=1}^\infty}
\newcommand{\pnn}{\prod_{n=1}^\infty}
\newcommand{\aln}{\alpha_n(\mu)}
\newcommand{\emp}{e^{\pi|Im\mu|}}
\newcommand{\ppp}{\prod_{p=p_0}^\infty}
\newcommand{\ppk}{\sum_{p=p_0}^{\infty}\sum_{k=1}^{[\ln p]}}
\newcommand{\pps}{\sum_{p=p_0}^\infty}
\newcommand{\mln}{m(\lambda_n)}
\newcommand{\mnk}{\mu_{n_k}}
\newcommand{\lnk}{\lambda_{n_k}}
\newcommand{\lqr}{\langle q \rangle}
\newcommand{\mlnk}{m(\lambda_{n_k})}

\newcommand{\sxm}{s(x,\mu)}
\newcommand{\skm}{s(\xi,\mu)}
\newcommand{\cxm}{c(x,\mu)}
\newcommand{\ckm}{c(\xi,\mu)}

\centerline {\bf On Properties of the Sturm-Liouville Operator with}
\centerline {\bf Degenerate Boundary Conditions}
\medskip
\medskip
\centerline {\bf Alexander Makin}
\medskip
\medskip

\medskip
\medskip
\begin{abstract}
We consider spectral problems for the Sturm-Liouville operator with arbitrary complex-valued potential $q(x)$ and degenerate boundary conditions. We solve corresponding inverse problem, and also study the completeness property and the basis property of the root function system.
\end{abstract}

{\bf 1. Introduction.}
Consider the Sturm-Liouville equation
$$
u''-q(x)u+\lambda u=0 \eqno (1)
$$
with two-point boundary conditions
$$
B_i(u)=a_{i1}u'(0)+a_{i2}u'(\pi)+a_{i3}u(0)+a_{i4}u(\pi)=0, \eqno (2)
$$
where the $B_i(u)$ $(i=1,2)$ are linearly independent forms with arbitrary
complex-valued coefficients and $q(x)$ is an arbitrary complex-valued
function of class $L_1(0,\pi)$.
It is convenient to write conditions (2) in the matrix form
$$
A=\left(
\begin{array}{cccc}
a_{11}&a_{12}&a_{13}&a_{14}\\
a_{21}&a_{22}&a_{23}&a_{24}
\end{array}
\right)
$$
and denote the matrix composed of the ith and jth columns of $A$
  $(1\le i<j\le4)$
by $A(ij)$; we set $A_{ij}=det A(ij)$.

It is known that conditions (2) can be divided into two classes:

1) nondegenerate conditions;

2) degenerate conditions.

Boundary conditions (2) are called  nondegenerate if they satisfy one of the
following relations:
$$
1) A_{12}\ne0,\quad2) A_{12}=0, A_{14}+A_{23}\ne0,\quad3) A_{12}=0,
 A_{14}+A_{23}=0, A_{34}\ne0.
$$
There is an enormous
literature related to the spectral theory of the Sturm-Liouville operator with nondegenerate boundary conditions.
In particular, the following assertion has been proved.

{\bf Theorem } ([1]). {\it For any
 nondegenerate conditions the
spectrum of problem (1), (2) consists of a countable set $\{\lambda_n\}$
of eigenvalues with only one limit point $\infty$, and the dimensions of the
corresponding root subspaces are bounded by one constant. The system $\{u_n(x)\}$
of eigen- and associated functions is complete and minimal in $L_2(0,1)$; hence,
it has a biorthogonally dual system $\{v_n(x)\}$.}

In this paper, we study eigenvalue problems for the Sturm-Liouville operator with degenerate boundary conditions.
This type of boundary conditions has been investigated much less.

{\bf 2. Preliminaries.}

Let boundary conditions (2) be degenerate. According to
[1, 2], this is equivalent to the fulfillment of the following conditions:

$$
A_{12}=0, \quad A_{14}+A_{23}=0,\quad A_{34}=0.
$$

According to [2], any boundary
conditions of the considered class
are equivalent to the boundary conditions determined by the matrix
$$
A=\left(
\begin{array}{cccc}
1&b&0&0\\
0&0&1&-b
\end{array}
\right),\quad
\mbox{or}\quad	A=\left(
\begin{array}{cccc}
0&1&0&0\\
0&0&0&1
\end{array}
\right).
$$
If in the first case $b=0$ then for any potential $q(x)$ we have the
initial value problem (the Cauchy problem) which has no eigenvalues. The same
situation takes place in the second case.

Further we will consider the first case if $b=\pm1$. Then the boundary conditions can be written
in more visual form
$$
u'(0)+(-1)^\theta u'(\pi)=0,\quad u(0)+(-1)^{\theta+1}u(\pi)=0. \eqno(3)
$$
$\theta=0,1$.
It is easily shown that if
$q(x)\equiv0$ then any $
 \lambda\in\mathbb{C}$ is an eigenvalue of infinite multiplicity. This abnormal
example illustrates the difficulty of investigation of problems with boundary conditions of the
considered class.

Denote by $c(x,\mu), s(x,\mu)$ $(\lambda=\mu^2)$ the fundamental system of solutions to
(1) with the initial conditions $c(0,\mu)=s'(0,\mu)=1$, $c'(0,\mu)=s(0,\mu)=0$.
The following identity is well known
$$
c(x,\mu)s'(x,\mu)-c'(x,\mu)s(x,\mu)=1.\eqno(4)
$$
Simple computations show that the characteristic equation of problem
(1), (3) can be reduced to the form $\Delta(\mu)=0$,
where
$$
\Delta(\mu)=c(\pi,\mu)-s'(\pi,\mu). \eqno(5)
$$
By $\Gamma(z,r)$ we denote the disk of radius $r$ centered at a point $z$.
By $PW_\sigma$ we denote the class of entire functions $f(z)$ of exponential type
$\le\sigma$ such that $||f(z)||_{L_2(R)}<\infty$, and by $PW_\sigma^-$
we denote the set of odd functions in $PW_\sigma$.

{\bf 3. Inverse problem.}

The following two assertions provide necessary and sufficient conditions to be
satisfied by the characteristic determinant $\Delta(\mu)$.

{\bf Theorem 1.}[3] {\it If a function $\Delta(\mu)$ is the characteristic determinant
of a problem (1), (3), then
$$
\Delta(\mu)=\frac{f(\mu)}{\mu},
$$
where $f(\mu)\in PW_\pi^-$.}

{\bf Theorem 2.} {\it Let a function $v(\mu)$ have the form
$$
v(\mu)=\frac{f(\mu)}{\mu},\eqno(6)
$$
where $f(\mu)\in PW_\pi^-$, and satisfies the condition
$$
\int_{-\infty}^{\infty}|\mu^{m}f(\mu)|^2d\mu<\infty, \eqno(7)
$$ where
$m$ is a nonnegative integer number.
Then, there exists a function $q(x)\in W_2^m(0,\pi)$ such that the characteristic
determinant of problem (1), (3) with the potential $q(x)$ satisfies $\Delta(\mu)=v(\mu)$.}

{\bf Proof.} If $m=0$ the theorem was proved in [3]. Further we will count that $m>0$. Since [4]
$$
|f(\mu)|\le C_1||f(\mu)||_{L_2(R)}e^{\pi|Im\mu|},\eqno(8)
$$
it follows that there exists an arbitrary large positive integer
$N$ such that
$$
|u(\mu)|<1/10, \quad |f(\mu)|<1\eqno(9)
$$
on the set $|Im\mu|\le1$, $Re\mu\ge N$.
Let $\mu_n$ $(n=1,2,\ldots)$ be a strictly monotone increasing sequence of positive
numbers such that $|\mu_n-(N+1/2)|<1/10$ if $1\le n\le N$
and $\mu_n=n$ if $n\ge N+1$. Consider the function
$$
s(\mu)=\pi\prod_{n=1}^\infty\frac{\mu_n^2-\mu^2}{n^2}=\spm\prod_{n=1}^N\frac{\mu_n^2-\mu^2}{n^2-\mu^2}.\eqno(10)
$$
Obviously, all zeros of the function $s(\mu)$ are simple, and, in addition, the
inequality
$$
(-1)^n\dot s(\mu_n)>0\eqno(11)
$$
holds for any
$n$.
Denote $P(\lambda)=\prod_{n=1}^N(\mu_n^2-\lambda)$, $Q(\lambda)=\prod_{n=1}^N(n^2-\lambda)$. Evidently,
$$
\frac{\lambda^m P(\lambda)}{Q(\lambda)}=\lambda^m+\sum_{j=1}^{m}\alpha_j\lambda^{m-j}+\frac{R(\lambda)}{Q(\lambda)},\eqno(12)
$$
where $R(\lambda)$ is a polynomial of degree  $N-1$ and $\alpha_j$ are some constants.
It follows from (12) that
$$
\prod_{n=1}^N\frac{\mu_n^2-\mu^2}{n^2-\mu^2}=1+\sum_{j=1}^{l}\alpha_j\mu^{-2j}+\frac{R(\mu^2)}{\mu^{2l}Q(\mu^2)}.\eqno(13)
$$
It follows from (10), (13) that
$$
\dot s(n)=\frac{\pi(-1)^n}{n}(1+\sum_{j=1}^{m}\alpha_j n^{-2j}+O(n^{-2m-2})).\eqno(14)
$$

Consider the equation
$$
z^2-u(\mu_n)z-1=0.\eqno(15)
$$
It has the roots
$$
c_n^\pm=\frac{u(\mu_n)\pm\sqrt{u^2(\mu_n)+4}}{2}.\eqno(16)
$$
It follows from (9) that for any $n$ all numbers $c_n^+$ lie in the disk $\Gamma(1,1/2)$ and
all numbers
$c_n^-$ lie in the disk $\Gamma(-1,1/2)$. Let for even
$n$ $c_n=c_n^+$, and for odd $n$ $c_n=c_n^-$.
Then $(-1)^nRec_n>0$
 for any $n=1, 2,\ldots$.
This, together with (11) implies that
$Rew_n>0$ for any $n$, where
$$
w_n=\frac{c_n}{\mu_n\dot s(\mu_n)}.\eqno(17)
$$

We set $F(x,t)=F_0(x,t)+\hat F(x,t)$, where
$$
F_0(x,t)=\sum_{n=1}^N\left(\frac{2c_n}{\mu_n\dot s(\mu_n)}
\sin\mu_nx\sin\mu_nt-
\frac{2}{\pi}\sin nx\sin nt\right),
$$
$$
\hat F(x,t)=\sum_{n=N+1}^\infty\left(\frac{2c_n}{\mu_n\dot s(\mu_n)}
\sin\mu_nx\sin\mu_nt-
\frac{2}{\pi}\sin nx\sin nt\right).\eqno(18)
$$
One can readily see that $F_0(x,t)\in C^\infty(R^2)$.
Consider the function $\hat F(x,t)$.
If $n\ge N+1$, then, by taking into account
(8), (16) and the rule for choosing the roots of equation (15), we obtain
$$
c_n=(-1)^n+\frac{f(n)}{2n}+(-1)^n(\sum_{j=1}^m \beta_j\frac{f^{2j}(n)}{n^{2j}}+O(1/n^{2m+2})).
\eqno(19)
$$
It follows from (8), (14), (18), and (19) that
$$
\begin{array}{c}
\hat F(x,t)=\sum_{n=N+1}^\infty\frac{2}{\pi}\left(\frac{1+(-1)^n\frac{f(n)}{2n}+\sum_{j=1}^m \beta_j\frac{f^{2j}(n)}{n^{2j}}+O(1/n^{2m+2})}{1+\sum_{j=1}^{l}\alpha_j n^{-2j}+O(n^{-2m-2})}-1
\right)\sin nx\sin nt=\\
=\sum_{n=N+1}^\infty\frac{2}{\pi}[(1+(-1)^n\frac{f(n)}{2n}+\sum_{j=1}^m\beta_j\frac{f^{2j}(n)}{n^{2j}}+O(1/n^{2m+2}) )\times\\(1+\sum_{j=1}^{m}\tilde\alpha_j n^{-2j}+O(n^{-2m-2}))-1
]\sin nx\sin nt=\\
=\frac{2}{\pi}\sum_{n=N+1}^\infty(\gamma_nf(n)+\sum_{j=1}^{m}\tilde\alpha_j n^{-2j}+O(n^{-2m-2}))\sin nx\sin nt=\\
=(\hat G(x-t)-\hat G(x+t))/2,

\end{array}
$$
where
$$
\hat G(y)=\frac{2}{\pi}\sum_{i=1}^3G_i(y), \quad G_1(y)=\sum_{n=N+1}^\infty \gamma_nf(n)\cos ny,
$$
$$
 G_2(y)=\sum_{n=N+1}^\infty\sum_{j=1}^{m}\tilde\alpha_j n^{-2j}\cos ny=
 \sum_{j=1}^{m}\tilde\alpha_j\sum_{n=N+1}^\infty n^{-2j}\cos ny,
 $$
 $$
 G_3(y)=\sum_{n=N+1}^\infty \tilde\gamma_n n^{-2m-2}\cos ny
$$

where $|\gamma_n|<C_2/n$, $|\tilde\gamma_n|<C_2$.

The relation
$$
\sum_{n=1}^\infty|nf(n)|^{2m}=\frac{1}{2}||\mu^mf(\mu)||_{L_2(R)},
$$ which follows from the Paley-Wiener theorem, together with the Parseval equality and (7), implies
that $G_1(y)\in W_2^{m+1}[0,2\pi]$.
It is known [5] that for any $j=1,2,\ldots$ the series $\sum_{n=N+1}^\infty n^{-2j}\cos ny$ are infinitely differentiable functions on the segment $[0, 2\pi]$. One can readily see that $G_3(y)\in W_2^{m+1}[0,2\pi]$.
 Therefore, we obtain the representation
$$
F(x,t)=F_0(x,t)+(\hat G(x-t)-\hat G(x+t))/2,\eqno(20)
$$
where the functions $F_0(x,t)$ and $\hat G(y)$ belong to the above-mentioned classes.

Now let us consider the Gelfand-Levitan equation
$$
K(x,t)+F(x,t)+\int_0^xK(x,s)F(s,t)ds=0\eqno(21)
$$
and prove that it has a unique solution in the space $L_2(0,x)$ for each $x\in[0,\pi]$.
To this end, it suffices  to show that the corresponding homogeneous equation
has only the trivial solution.

Let $f(t)\in L_2(0,x)$. Consider the equation
$$
f(t)+\int_0^xF(s,t)f(s)ds=0.
$$
Following [6], by multiplying the last equation by $\bar f(t)$ and by integrating the
resulting relation over the interval
$[0,x]$, we obtain
$$\begin{array}{c}
\int_0^x|f(t)|^2dt+
\sum_{n=1}^\infty\frac{2c_n}{\mu_n\dot s(\mu_n)}\int_0^x\bar f(t)\sin\mu_ntdt\int_0^xf(s)\sin\mu_nsds-\\
-\sum_{n=1}^\infty\frac{2}{\pi}\int_0^x\bar f(t)\sin ntdt\int_0^xf(s)\sin nsds=0.
\end{array}
$$
This, together with the Parseval equality for the function system $\{\sin nt\}_1^\infty$
on the interval $[0,\pi]$ implies that
$$
\sum_{n=1}^\infty w_n|\int_0^xf(t)\sin\mu_ntdt|^2=0,
$$
where the $w_n$ are the numbers given by (17). Since $Rew_n>0$, we see that
$\int_0^xf(t)\sin\mu_ntdt=0$ for any $n=1,2,\ldots$. Since [7, 8]
the system $\{\sin\mu_nt\}_1^\infty$ is complete on the interval $[0,\pi]$, we have $f(t)\equiv0$ on
$[0,x]$.

Let $\hat K(x,t)$ be a solution of equation (21), and let $\hat q(x)=2\frac{d}{dx}\hat K(x,x)$;
then it follows [9] from (20) that $\hat q(x)\in W_2^m(0,\pi)$. By $\hat s(x,\mu)$, $\hat c(x,\mu)$
we denote the fundamental solution system of equation (1) with potential $\hat q(x)$ and
the initial conditions $\hat s(0,\mu)=\hat c'(0,\mu)=0$, $\hat c(0,\mu)=\hat s'(0,\mu)=1$.
By reproducing the corresponding considerations in [6], we obtain $\hat s(\pi,\mu)\equiv s(\mu)$, whence
it follows that the numbers $\mu_n^2$ form the spectrum of the Dirichlet problem for equation (1) with
potential $\hat q(x)$, and $\hat c(\pi,\mu_n)=c_n$, which, together with
identity (4), implies that
$\hat s'(\pi,\mu_n)=1/c_n$.

Let $\hat\Delta(\mu)$ be the characteristic determinant of problem (1), (3) with potential $\hat q(x))$. Let us prove that $\hat\Delta(\mu)\equiv v(\mu)$.
By theorem 1, the function $\hat \Delta(\mu)$ admits the representation
$$
\hat\Delta(\mu)=\frac{\hat f(\mu)}{\mu},
$$
where $\hat f(\mu)\in PW_\pi^-$.
By taking into account relation (4) and the fact that the numbers $c_n$ are roots of equation (15),
we have
$$
\hat\Delta(\mu_n)=
\hat c(\pi,\mu_n)-\hat s'(\pi,\mu_n)
=c_n-c_n^{-1}=v(\mu_n).
$$
It follows that the function
$$
\Phi(\mu)=\frac{u(\mu)-\hat\Delta(\mu)}{s(\mu)}=\frac{f(\mu)-\hat f(\mu)}{\mu s(\mu)}
$$
is an entire function on the complex plane. Since the function $g(\mu)=f(\mu)-\hat f(\mu)$
belongs to $\pwp$,
it follows from (8) that
$$
|g(\mu)|\le C_3e^{\pi|Im\mu|}.\eqno(22)
$$
From (10), we find that if $|Im\mu|\ge1$, then
$$
|\mu s(\mu)|\ge C_4e^{\pi|Im\mu|}\eqno(23)
$$
$(C_4>0)$.If $|Im\mu|\ge1$, then we obtain the estimate $|\Phi(\mu)|\le C_3/C_4$.

By $H$ we denote the union of the vertical segments $\{z: |Rez|=n+1/2, |Imz|\le1\}$,
where $n=N+1,N+2,\ldots$. It follows from (10) that if $\mu\in H$, then
$|\mu s(\mu)|\ge C_5>0$. The last inequality, together with (22), (23), and the
maximum principle for the absolute value of an analytic function, implies that
 $|\Phi(\mu)|\le C_6$ in the strip $|Im\mu|\le1$.
Consequently, the function $\Phi(\mu)$ is bounded on the entire complex plane;
therefore, by the Liouville theorem, it is a constant . It follows from the Paley-Wiener theorem and the
Riemann lemma [1] that if $|Im\mu|=1$, then $\lim_{|\mu|\to\infty}g(\mu)=0$,
whence, we obtain $\Phi(\mu)\equiv0$.

{\bf 4. Completeness and the basis property.}

Completeness of the root function system of problem (1), (3) was investigated in [10]. In particular, it was shown that if $q(x) \in C^k[0,\pi]$ for some $k\ge0$,
and $q^{(k)}(0)\ne(-1)q^{(k)}(\pi)$, then the root function system is complete in $L_2(0,\pi)$. If
there exists an $\varepsilon>0$ such that $q(x)-q(\pi-x)=0$
for almost all $x\in[0,\varepsilon]$, then the mentioned system is not complete in $L_2(0,\pi)$.
In was established in [3], [11] that there exist potentials $q(x)$ such that the root function systems of corresponding problems (1), (3) are complete in $L_2(0,\pi)$ and contain associated functions of arbitrary high order, i.e. the dimensions of root subspaces infinitely grow.

Since for a wide class of potentials $q(x)$ the root function system of problem (1), (3) is complete in $L_2(0,\pi)$ one can set a question whether the mentioned system forms a basis.

Let $\lambda_n=\mu_n^2$ $(Re\mu_n\ge0, n=1,2,\ldots)$ be the eigenvalues of problem (1), (3) numbered neglecting their multiplicities in nondecreasing order of absolute value. By $\mln$ we denote the multiplicity of an eigenvalue $\lambda_n$. In addition, assume that
the function $q(x)$ is continuous on the interval $(0,\pi)$.

{\bf Theorem 3.} {\it Suppose a subsequence of eigenvalues $\lnk$  satisfies  the following two conditions:

1. $|Im\mnk|<M$;

2. $\lim_{k\to\infty}\frac{m(\lnk)}{\ln|\lnk|}=0$;

 Then the system of eigenfunctions and associated functions of problem (1), (3) is not a basis in $L_2(0,\pi)$.}

{\bf Proof.} Let us calculate the Green function $G(x,\xi,\mu)$ of operator (1), (3). By [9], $G(x,\xi,\mu)=H(x,\xi,\mu)/\dm$, where
$H(x,\xi,\mu)=\Phi(x,\xi,\mu)/2+g(x,\xi)\dm$, where
$$
\begin{array}{c}															
\Phi(x,\xi,\mu)=\sxm\{c'(\pi,\mu)[-\ckm s(\pi,\mu)-\skm(-1-c(\pi,\mu))]-\\-
[1-c(\pi,\mu)][\ckm(-1+s'(\pi,\mu))-\skm c'(\pi,\mu)]\}-\\								
-\cxm\{[1+s'(\pi,\mu)][-\ckm s(\pi,\mu)-\skm(-1-c(\pi,\mu))]+\\+			
s(\pi,\mu)[\ckm(-1+s'(\pi,\mu))-\skm c'(\pi,\mu)]\}, 																						 
\end{array}\eqno(24)																	
$$																																																																																																																												 
$g(x,\xi)=\pm(s(x,\mu)c(\xi,\mu)-c(x,\mu)s(\xi,\mu))/2$,
 the sigh $"+"$ is used for $x>\xi$, and the sigh $"$ $-$ $"$ is used for $x<\xi$. Combining like terms in (24) gives
$$
\begin{array}{c}
\Phi(x,\xi,\mu)=2[\sxm\ckm-\cxm\skm]-\\-
[c(\pi,\mu)+s'(\pi,\mu)][\sxm\ckm+\cxm\skm]+\\+2[c'(\pi,\mu)\sxm\skm+s(\pi,\mu)\cxm\ckm].
\end{array}\eqno(25)
$$

 Let $e(x,\mu)$ be the solution of equation (1) satisfying the initial conditions $e(0,\mu)=1$, $e'(0,\mu)=i\mu$, and let $K(x,t)$, $K^+(x,t)=K(x,t)+K(x,-t)$, and
$K^-(x,t)=K(x,t)-K(x,-t)$ be the transformation kernels [1] realizing the representations
$$
e(x,\mu)=e^{i\mu x}+\int_{-x}^xK(x,t)e^{i\mu t}dt,
$$
$$
c(x,\mu)=\cos\mu x+\int_0^xK^+(x,t)\cos\mu tdt,
$$
$$
s(x,\mu)=\frac{\sin\mu x}{\mu}+\int_0^xK^-(x,t)\frac{\sin\mu t}{\mu}dt.\eqno(26)
$$
It was shown in [12] that
$$
c(\pi,\mu)=\cos\pi\mu+\frac{\pi}{2}\lqr\spm-\int_0^\pi\frac{\partial K^+(\pi,t)}{\partial t}\frac{\sin\mu t}{\mu}dt,\eqno(27)
$$

$$
s'(\pi,\mu)=\cos\pi\mu+\frac{\pi}{2}\lqr\spm+\int_0^\pi\frac{\partial K^-(\pi,t)}{\partial x}\frac{\sin\mu t}{\mu}dt,\eqno(28)
$$
where $\lqr=\frac{1}{\pi}\int_0^\pi q(x)dx$.
By differentiating the second of equalities (26) and taking into account [12] that $K^+(\pi,\pi)=\frac{\pi}{2}\lqr$, we obtain
$$
c'(\pi,\mu)=-\mu\sin\pi\mu+
\frac{\pi}{2}\lqr\cos\pi\mu+\int_0^\pi\frac{\partial K^+(\pi,t)}{\partial x}\cos\mu tdt.\eqno(29)
$$
By substituting the right-hand sides of (26-29) in (25), we get
$$
\begin{array}{c}
\Phi(x,\xi,\mu)=2(\sin\mu x\cos\mu\xi-\cos\mu x\sin\mu\xi)/\mu-\\-
2\cos\pi\mu(\sin\mu x\cos\mu\xi+\cos\mu x\sin\mu\xi)/\mu+\\+
2(-\sin\pi\mu\sin\mu x\sin\mu\xi+\sin\pi\mu\cos\mu x\cos\mu\xi)/\mu+o(\mu^{-1})\emp=\\=
2[\sin\mu(x-\xi)+\sin\mu(\pi-(x+\xi))]/\mu+o(\mu^{-1})\emp.
\end{array}
$$

Throughout the following we assume that $|Im\mu|<M$. Then the last equality implies that
$$ G(x,\xi,\mu)=\frac{R(x,\xi,\mu)}{\dm}+g(x,\xi),
\eqno(30)
$$
where
$$
R(x,\xi,\mu)=[\sin\mu(x-\xi)+\sin\mu(\pi-(x+\xi))]/\mu+o(\mu^{-1}).\eqno(31)
$$

Let us study the function $G(x,\xi,\mu)$ in the  neighborhood of the eigenvalues $\lambda_n$. It follows from [3]
that each root subspace contains one eigenfunction and possibly associated functions. Let
$\{\stackrel{h}{u}_n(x)\}$ $(h=\overline{0,m(\lambda_n)})$ be an arbitrary canonical system of eigenfunctions and associated functions of problem (1), (3), and let $\{\stackrel{h}{v}_n(x)\}$  be appropriately normalized canonical system of eigenfunctions and associated functions of the adjoint boundary value problem [13],
 i.e. $\stackrel{0}{u}_n(x)$ and $\stackrel{0}{v}_n(x)$ are eigenfunctions, and
$\stackrel{h}{u}_n(x)$ and $\stackrel{h}{v}_n(x)$ $(h\ge1)$ are associated functions of order $h$, where	
$$
(\stackrel{h}{u}_n(x),\stackrel{g}{v}_k(x))_{L_2(0,\pi)}=\delta_{n,k}\delta_{h,m(\lambda_n)-1-g}.
$$
Further we consider only root subspaces corresponding the mentioned-above subsequence of the eigenvalues $\lambda_{n_k}$.
Denote
$$
\stackrel{0}{R}_{n_k}(x,\xi)=\stackrel{0}{u}_{n_k}(x)\overline{\stackrel{0}{v}_{n_k}(\xi)},
$$
$$
\stackrel{\mlnk-1}{R}_{n_k}(x,\xi)=\sum_{p=0}^{\mlnk-1}\stackrel{p}{u}_{n_k}(x)\overline{\stackrel{\mlnk-1-p}{v}_{n_k}(\xi)}.
$$

Since the function $f(\mu)$  has a root of multiplicity $\mlnk$ at the point $\mnk$, then
$$
f(\mu)=\sum_{l=\mlnk}^\infty c_l(\mu-\mnk)^l=(\mu-\mnk)^{\mlnk}\sum_{l=0}^\infty c_{\mlnk+l}(\mu-\mnk)^l.\eqno(32)
$$
Obviously, $c_{\mlnk}=f^{(\mlnk)}(\mnk)/\mlnk!$. Relations (30) and (32), together with [13] imply the equality
$$
\begin{array}{c}
\stackrel{0}{R}_{n_k}(x,\xi)=\lim_{\mu\to\mnk}(\mu^2-\mnk^2)^{\mlnk} G(x,\xi,\mu)=\\=
\frac{2^{\mlnk} \mlnk!\mnk^{\mlnk+1}R(x,\xi,\mnk)}{f^{(\mlnk)}(\mnk)}.
\end{array}\eqno(33)
$$
From the Bernstein inequality [4], we obtain
$$
|f^{(\mlnk)}(\mnk)|\le C_1\pi^{\mlnk}.\eqno(34)
$$
It follows from (33) and (34) that
$$
|\stackrel{0}{R}_{n_k}(x,\xi)|\ge C_2(2/\pi)^{\mlnk}\mlnk!|\mnk|^{\mlnk+1}|R(x,\xi,\mnk)|,
$$
where $C_2>0$, hence,
$$
\begin{array}{c}
||\stackrel{0}{R}_{n_k}(x,\xi))||^2_{L_2((0,\pi)\times(0,\pi))}\ge\\\ge C_3[(2/\pi)^{\mlnk}\mlnk!]^2|\mnk|^{2\mlnk+2}||R(x,\xi,\mnk)||^2_{L_2((0,\pi)\times(0,\pi))},
\end{array}\eqno(35)
$$
where $C_3>0$.
																		
If $\stackrel{0}{v}_{n_k}(\xi)\ne0$, then the function $\stackrel{\mlnk-1}{R}_{n_k}(x,\xi)$ is an associated  function of order $\mlnk-1$,
corresponding to the eigenvalue $\mnk^2$ and the eigenfunction $\stackrel{0}{u}_{n_k}(x)$. It follows from [14] that
$$
\begin{array}{c}
||\stackrel{0}{R}_{n_k}(x,\xi)||^2_{L_2(\pi/3,\pi/2)}\le \\ \le [C_4\mlnk|\mnk|]^{2\mlnk-2}||\stackrel{\mlnk-1}{R}_{n_k}(x,\xi)||^2_{L_2(\pi/4,3\pi/4)},
\end{array}
\eqno(36)
$$
where $C_4$ is a constant independent of $\xi$. If $\stackrel{0}{v}_{n_k}(\xi)=0$,	then the validity of (36) is obvious.
It follows from (36) and [15] that

$$
\begin{array}{c}
||\stackrel{0}{R}_{n_k}(x,\xi)||^2_{L_2(0,\pi)}\le
[C_5\mlnk|\mnk|]^{2\mlnk-2}
||\stackrel{\mlnk-1}{R}_{n_k}(x,\xi)||^2_{L_2(0,\pi)}.
\end{array}\eqno(37)
$$

 By integrating inequality (37) with respect to $\xi$, we have
$$
\begin{array}{c}
||\stackrel{0}{R}_{n_k}(x,\xi)||^2_{L_2((0,\pi)\times(0,\pi)}\le\\\le
[C_6\mlnk|\mnk|]^{2\mlnk-2}
||\stackrel{\mlnk-1}{R}_{n_k}(x,\xi)||^2_{L_2((0,\pi)\times(0,\pi))}.
\end{array}\eqno(38)
$$
It follows from (31) that
$$
||R(x,\xi,\mnk)||^2_{L_2((0,\pi)\times(0,\pi))}\ge C_7|\mnk|^{-2},\eqno(39)
$$
where $C_7>0$.

Relations (35), (38), (39), and the Stirling formula imply that
$$
||\stackrel{\mlnk-1}{R}_{n_k}(x,\xi)||^2_{L_2((0,\pi)\times(0,\pi))}\ge \frac{(\mlnk!|\mnk|)^2}{(C_8^{\mlnk}\mlnk)^{2\mlnk}}\ge\frac{|\mnk|^2}{C_9^{\mlnk}} ,
$$
where $C_8, C_9>0$. By the conditions of the theorem, the right-hand side of the last inequality tends to infinity as $k\to\infty$. This, combined with a resonance type theorem
[16] implies the validity of theorem 3.

\medskip
\medskip

\centerline {Acknowledgements}
\medskip

The research was supported by the Russian Foundation for Basic Research (project no. 10-01-00411).
\medskip
\medskip

\centerline {\bf References}
\medskip
\medskip
\begin{itemize}
\item[1.] V.A. Marchenko. Sturm-Liouville Operators and Their Applications. Kiev, 1977
(in Russian); English transl.: Birkh\"{a}user, Basel, 1986.

\item[2.] P. Lang and J. Locker. Spectral theory of two-point differential operators
determined by $-D^2$. II. Analysis of cases. J. Math. Anal. Appl. {\bf 146} (1990), 148-191.

\item[3.] Makin A.S. Characterzation of the spectrum of the Sturm-liouville operator with irregular boundary conditions. Differ.Equations, 46, No. 10, 1427-1437 (2010); translation from Differ. Uravn., 46, No. 10, 1421-1432 (Russian) (2010).

\item[4.] S.M. Nikolskii. Approximation of Functions of Several Variables and Embedding Theorems
(Moscow, 1977).

\item[5.] A.P. Prudnikov, U.A. Brychkov, O.I. Marichev. Integrals and series. Moscow, Nauka, 1981.

\item[6.] V.A. Tkachenko. Discriminants and generic spectra of nonselfadjoint Hill's operators,
 Adv. Sov. Mathem. {\bf 19}, 41-71 (1994).

\item[7.] A.M. Sedletskii. The stability of the completeness and of the minimality in $L_2$ of a
system of exponential functions, Mat. Zametki {\bf 15}, 213-219 (1974).

\item[8.] A.M. Sedletskii. Convergence of angarmonic Fourier series in systems of exponentials,
cosines and sines, Dokl. Akad. Nauk SSSR {\bf 301}, 1053-1056 (1988).

\item[9.] [8] M.A. Naimark. Linear Differential Operators (in Russian). Nauka, Moscow, 1969; English transl.:
 Ungar, New-York, 1967.

\item[10.] M.M. Malamud. On the Completeness of the System of Root Vectors of the Sturm--Liouville Operator
with General Boundary Conditions, Funct. An. and Its Appl. V. 42, No. 3, 198-204 (2008); translation from
Funkt. An. i ego Pr. 42, No. 3, 45-52 (Russian) (2008).

\item[11.] A.S. Makin, Differ. Equat. On a two-point boundary value problem for the Sturm-Liouville
operator with nonclassical spectrum asymptotics, accepted.

\item[12.] L.A. Pastur, V.A. Tkachenko. Spectral theory of a class of one-dimensional
Schrodinger operators with limit-peroodic potentials, Trudy Moscow Mat. Obshch.
{\bf 51}, 114-168 (1988).

\item[13.] M. V. Keldysh. On the Completeness of eigenfunctions of some classes of
nonselfadjoint linear operators // Uspekxi Mat. Nauk 1971. V. 26, No. 4., 15-41.

\item[14.]  A.S. Makin. Exact estimates for the root functions of a second order general elliptic operator.
Doklady Mathematics, 57, No. 1, 98-100 (1998); translation from Dokl. Akad. Nauk, 358, No. 5, 594-596 (Russian)
 (1998).

\item[15.] I.S. Lomov. Some properties of eigen- and associated functions of the
Sturm--Liouville operator, Differ. Equations 18, No. 10, 1206--1213 (1982);
translation from Differ. Uravn. 18, No. 10, 1684--1694, (Russian) (1982).

\item[16.] K. Iosida. Functional analysis. Moscow, Mir, 1967 (in Russian).

\end{itemize}
\medskip
\medskip
Moscow State University of Instrument-Making and Computer Science,
Stromynka 20, Moscow, 107996, Russia
\medskip

\medskip

E-mail address: alexmakin@yandex.ru

}
\end{document}